\documentclass[a4paper,12pt]{article}
\usepackage[dvips]{epsfig,graphicx}
\usepackage{amsmath,amsfonts,verbatim,latexsym}

\newcommand{\RR}{{\mathbb R}}

\newcommand{\qed}{\hfill\hbox{\rule{3pt}{6pt}}}
\newcommand{\proof}{{\sc Proof. }}
\newcommand{\G}{\Gamma}
\newcommand{\mat}{{\rm Mat}}

\newcommand{\as}{\alpha^*}
\renewcommand{\a}{\alpha}

\newcommand{\la}{\langle}
\newcommand{\ra}{\rangle}

\newcommand{\MX}{\mat_X(\RR)}
\newcommand{\Mdva}{\mat_2(\RR)}

\newcommand{\Al}{{\cal L}}
\newcommand{\Sal}{{\cal L}^{sym}}
\newcommand{\Asal}{{\cal L}^{asym}}

\newtheorem{theorem}{Theorem}[section]
\newtheorem{lemma}[theorem]{Lemma}
\newtheorem{corollary}[theorem]{Corollary}
\newtheorem{proposition}[theorem]{Proposition}
\newtheorem{definition}[theorem]{Definition}

\newtheorem{example}[theorem]{Example}

\title{The $A$-like matrices for a hypercube}

\author{\v{S}tefko Miklavi\v{c} \\ 
        UP PINT and UP FAMNIT \\
        University of Primorska \\
        Muzejski trg 2 \\ 
        6000 Koper, Slovenia \\
        stefko.miklavic@upr.si \and
        Paul Terwilliger \\
        Department of Mathematics \\
        University of Wisconsin \\
        480 Lincoln Drive \\
        Madison WI 53706-1388, USA \\
        terwilli@math.wisc.edu}

\begin{document}
\maketitle

\begin{abstract}
  Let $D$ denote a positive integer and let $Q_D$ denote the graph of the $D$-dimensional
  hypercube. Let $X$ denote the vertex set of $Q_D$ and let $A \in \MX$ denote the adjacency
  matrix of $Q_D$. A matrix $B \in \MX$ is called $A$-{\em like} whenever both
  (i) $BA = AB$; (ii) for all $x,y \in X$ that are not equal or adjacent, the $(x,y)$-entry of $B$ is zero.
  Let $\Al$ denote the subspace of $\MX$ consisting of the $A$-like elements.
  We decompose $\Al$ into the direct sum of its symmetric part and antisymmetric part.
  We give a basis for each part. The dimensions of the symmetric part and antisymmetric part
  are $D+1$ and ${D \choose 2}$, respectively.
\end{abstract}

\section{Introduction}
\label{sec:intro}

Let $\G$ denote a finite undirected graph with vertex set $X$. 
Let $A \in \MX$ denote the adjacency matrix of $\G$.
For $B \in \MX$, we define $B$ to be $A$-{\em like} whenever both
\begin{itemize}
\item[(i)]
$BA = AB$;
\item[(ii)]
for all $x,y \in X$ that are not equal or adjacent, the $(x,y)$-entry of $B$ is zero.
\end{itemize}
Let $\Al$ denote the subspace of $\MX$ consisting of the $A$-like elements for $\G$.
We will discuss $\Al$ after a few comments.

\noindent
For $B \in \MX$ let $B^t$ denote the transpose of $B$.
Recall that $B$ is {\em symmetric} (resp. {\em antisymmetric}) whenever
$B^t=B$ (resp. $B^t=-B$).
For a subspace $H \subseteq \MX$, the set of symmetric (resp. antisymmetric)
matrices in $H$ is a subspace of $H$ called the {\em symmetric part} 
(resp. {\em antisymmetric part}) of $H$.
Note that the following are equivalent: 
(i) $H$ is the direct sum of its symmetric and antisymmetric part;
(ii) $H$ is closed under the transpose map.

\noindent
Let $\Sal$ (resp. $\Asal$) denote the symmetric (resp. antisymmetric) part of $\Al$.
Note that $\Al$ is closed under the transpose map, so the sum $\Al = \Sal + \Asal$
is direct.

\noindent
For a positive integer $D$ let $Q_D$ denote the graph of 
the $D$-dimensional hypercube (see Section \ref{sec:hyper} for formal definitions).  
For $\G=Q_D$ we find a basis for $\Sal$ and $\Asal$.
In particular we show that the dimensions of $\Sal$ and $\Asal$ are $D+1$ and ${D \choose 2}$, respectively.


\section{Preliminaries}
\label{sec:prelim}

Let $X$ denote a nonempty finite set. Let $\MX$ denote the $\RR$-algebra consisting of the matrices 
with entries in $\RR$, and rows and columns indexed by $X$. Let $V=\RR^X$ denote the vector space over 
$\RR$ consisting of the column vectors with entries in $\RR$ and rows indexed by $X$. Observe that $\MX$ 
acts on $V$ by left multiplication. We refer to $V$ as the {\em standard module} of $\MX$. For 
$v \in V$ let $v^t$ denote the transpose of $v$.
We endow $V$ with the bilinear form $\la u,v \ra = u^t v \; (u,v \in V)$. 
The vector space $V$ together with $\la \, , \, \ra$ is a Euclidean space.
For $x \in X$ let $\hat{x}$ denote the vector in $V$ that has $x$-coordinate $1$ and all other 
coordinates $0$. Observe that $\{\hat{x} | x \in X\}$ is an orthonormal basis for $V$.
For $B \in \MX$ we have $\la Bu,v \ra = \la u, B^t v \ra$ for all $u,v \in V$.
Therefore $B$ is symmetric (resp. antisymmetric) if and only if $\la Bu, v \ra = \la u, Bv \ra$
(resp. $\la Bu, v \ra = - \la u, Bv \ra$) for all $u,v \in V$.

\smallskip \noindent
Given a subspace $U \subseteq V$, let ${\rm End}(U)$ denote the $\RR$-algebra consisting of the linear 
transformations from $U$ to $U$. We identify ${\rm End}(V)$ with $\MX$.

\begin{lemma}
\label{lema_1}
For a subspace $U \subseteq V$ and $B \in {\rm End}(U)$ there exists a
unique $B^\dagger \in {\rm End}(U)$ such that $\la Bu, v \ra = \la u, B^\dagger v \ra$
for all $u,v \in U$. We call $B^\dagger$ the {\em adjoint} of $B$ relative to $\la \, , \, \ra$.
\end{lemma}
\proof
By the Fischer-Riesz theorem \cite[Theorem 9.18]{Rom}, for all $v \in U$ there exists a unique $z_v \in U$ such that
$\la Bu, v \ra = \la u, z_v \ra$ for all $u \in U$. 
Observe that the map $U \to U$, $v \mapsto z_v$ is an element of ${\rm End}(U)$
which we denote by $B^\dagger$. By construction $\la B u, v \ra = \la u, B^\dagger v \ra$ for all $u,v \in U$.
We have shown that $B^\dagger$ exists. 
The uniqueness of $B^\dagger$ follows from the uniqueness of the vector $z_v$ in our preliminary remarks. \qed

\begin{example}
Referring to Lemma \ref{lema_1}, assume $U=V$. Then $B^\dagger=B^t$.
\end{example}

\medskip \noindent
The lemma and corollaries below follow from Lemma \ref{lema_1}.

\begin{lemma}
\label{lema_2}
Let $U$ denote a subspace of $V$ and fix an orthonormal basis for $U$.
Pick $B \in {\rm End}(U)$ and consider the matrices representing $B$ and $B^\dagger$
with respect to this basis. Then these matrices are transposes of each other.
\end{lemma}

\begin{corollary}
\label{pos_3}
For a subspace $U \subseteq V$ and $B \in {\rm End}(U)$ the following 
{\rm (i)--(iv)} are equivalent.
\begin{itemize}
\item[{\rm (i)}]
There exists an orthonormal basis for $U$ with respect to which the matrix representing $B$ is symmetric.
\item[{\rm (ii)}]
With respect to any orthonormal basis for $U$ the matrix representing $B$ is symmetric.
\item[{\rm (iii)}]
$\la Bu, v \ra = \la u, Bv \ra$ for all $u,v \in U$.
\item[{\rm (iv)}]
$B=B^\dagger$.
\end{itemize}
\end{corollary}

\begin{corollary}
\label{pos_4}
For a subspace $U \subseteq V$ and $B \in {\rm End}(U)$ the following 
{\rm (i)--(iv)} are equivalent.
\begin{itemize}
\item[{\rm (i)}]
There exists an orthonormal basis for $U$ with respect to which the matrix representing $B$ is antisymmetric.
\item[{\rm (ii)}]
With respect to any orthonormal basis for $U$ the matrix representing $B$ is antisymmetric.
\item[{\rm (iii)}]
$\la Bu, v \ra = - \la u, Bv \ra$ for all $u,v \in U$.
\item[{\rm (iv)}]
$B = - B^\dagger$.
\end{itemize}
\end{corollary}
Motivated by Corollaries \ref{pos_3} and \ref{pos_4} we make a definition.

\begin{definition}
\label{def_1}
{\rm For a subspace $U \subseteq V$ and $B \in {\rm End}(U)$ we call $B$ {\it symmetric} (resp. {\it antisymmetric})
whenever the equivalent conditions {\rm (i)--(iv)} hold in Corollary \ref{pos_3} (resp. Corollary \ref{pos_4}).}
\end{definition}

\begin{definition}
\label{def_2}
{\rm Pick a subspace $U \subseteq V$ and a subspace $H \subseteq {\rm End}(U)$.
By the {\it symmetric part} of $H$ we mean the subspace of $H$ consisting of its symmetric elements.
By the {\it antisymmetric part} of $H$ we mean the subspace of $H$ consisting of its antisymmetric elements.}
\end{definition}
We make two observations.
\begin{lemma}
\label{lema_5}
For a subspace $U \subseteq V$ and a subspace $H \subseteq {\rm End}(U)$
the following {\rm (i), (ii)} are equivalent.
\begin{itemize}
\item[{\rm (i)}]
$H$ is the direct sum of its symmetric and antisymmetric part.
\item[{\rm (ii)}]
$H$ is closed under the adjoint map.
\end{itemize}
\end{lemma}
\begin{lemma}
\label{lem:prelim2}
Let $B \in \MX$ and assume $B$ is symmetric (resp. antisymmetric).
Then for any $B$-invariant subspace $U \subseteq V$, the restriction of 
$B$ to $U$ is symmetric (resp. antisymmetric) in the sense of Definition \ref{def_1}.
\end{lemma}

\noindent
Let $\G=(X,R)$ denote a finite, undirected, connected graph, without loops or multiple edges, with 
vertex set $X$, edge set $R$, path-length distance function $\partial$, and diameter 
$D:=\max \{\partial(x,y) |$ $\: x,y \in X\}$. For a vertex $x \in X$ and an integer $i \ge 0$ let
$\G_i(x) = \{ y \in X \mid \partial(x,y)=i \}$. 
We abbreviate $\G(x) = \G_1(x)$. For an integer $k \ge 0$ we say 
$\G$ is {\em regular with valency} $k$ whenever $|\G(x)|=k$ for all $x \in X$. We say $\G$ is 
{\em distance-regular} whenever for all integers $0 \le h,i,j \le D$ and all $x,y \in X$ with 
$\partial(x,y)=h$ the number $p_{ij}^h := | \G_i(x) \cap \G_j(y) |$
is independent of $x,y$. The constants $p_{ij}^h$ are known as the {\em intersection numbers} of $\G$. 
From now on assume $\G$ is distance-regular with $D \ge 1$. 
Observe that $\G$ is regular with valency $k=p_{11}^0$.

\noindent
We now recall the Bose-Mesner algebra of $\G$. For $0 \le i \le D$ let $A_i$ denote the 
matrix in $\MX$ with entries
$$
  (A_i)_{x y} = \left\{ \begin{array}{ll}
                 1 & \hbox{if } \; \partial(x,y)=i, \\
                 0 & \hbox{if } \; \partial(x,y) \ne i \end{array} \right. \qquad (x,y \in X).
$$
We abbreviate $A=A_1$ and call this the {\em adjacency matrix of} $\G$. Let $M$ 
denote the subalgebra of $\MX$ generated by $A$. By \cite[Theorem 11.2.2]{God} the matrices
$\lbrace A_i \rbrace_{i=0}^D$ form a basis for $M$. We call $M$ the 
{\em Bose-Mesner algebra of} $\G$. Observe that $M$ is commutative and semi-simple. 
By \cite[Theorem 12.2.1]{God} there exists a basis $\lbrace E_i \rbrace_{i=0}^D$ for $M$ such that
(i) $E_0 = |X|^{-1} J$,
(ii) $I = \sum_{i=0}^D E_i$,
(iii) $E_i^t = E_i \; (0 \le i \le D)$,
(iv) $E_i E_j = \delta_{ij} E_i \; (0 \le i,j \le D)$,
where $I$ and $J$ denote the identity and the all-ones matrix of $\MX$, respectively.
The matrices $\lbrace E_i \rbrace_{i=0}^D$ are known as the 
{\em primitive idempotents of} $\G$, and $E_0$ is called the {\em trivial} idempotent.
We recall the eigenvalues of $\G$. Since $\lbrace E_i \rbrace_{i=0}^D$ is a basis for $M$, there exist
real scalars $\lbrace \theta_i \rbrace_{i=0}^D$ such that
\begin{equation}
\label{theta}
  A = \sum_{i=0}^D \theta_i E_i.
\end{equation}
Combining this with (iv) above we obtain $AE_i = E_iA =\theta_i E_i$ for 
$0 \le i \le D$. Using (i) above we find $\theta_0=k$.
For $0 \le i \le D$ we call $\theta_i$ the {\em eigenvalue} of $\G$ associated with $E_i$.
The eigenvalues $\{\theta_i\}_{i=0}^D$ are mutually distinct since $A$ 
generates $M$. For $0 \le i \le D$ let $m_i$ denote the rank of $E_i$. 
We call $m_i$ the {\em multiplicity} of $\theta_i$. 

\smallskip \noindent
By (ii)--(iv) above, 
\begin{equation}
\label{de}
  V = E_0V + E_1V + \cdots + E_DV \qquad \hbox{(orthogonal direct sum)}.
\end{equation}
For $0 \le i \le D$ the space $E_iV$ is the eigenspace of $A$ associated with $\theta_i$.

\smallskip \noindent


\section{The tensor product}
\label{sec:kronecker}

In this section we recall the tensor product of vectors and matrices. 

\noindent
For $v \in \RR^X$ and $v' \in \RR^{X'}$ let $v \otimes v'$ denote the vector in 
$\RR^{X \times X'}$ with a $(x,x')$-entry equal to the $x$-entry of $v$ times the 
$x'$-entry of $v'$. We call $v \otimes v'$ the {\em tensor product} of $v$ and $v'$.

\noindent
For $B \in \MX$ and $B' \in \mat_{X'}(\RR)$ let $B \otimes B'$ denote the matrix in 
$\mat_{X \times X'}(\RR)$ with a $((x,x'),(y,y'))$-entry equal to the $(x,y)$-entry of $B$ times the 
$(x', y')$-entry of $B'$. We call $B \otimes B'$ the {\em tensor product} of $B$ and $B'$.
Pick $B_1, B_2 \in \MX$ and $B'_1, B'_2 \in \mat_{X'}(\RR)$.
Then by \cite[p.~107]{Eves},
\begin{equation}
\label{kron1}
  (B_1 \otimes B'_1)(B_2 \otimes B'_2) = (B_1 B_2) \otimes (B'_1 B'_2).
\end{equation}
Similarly, for $v_1, v_2 \in \RR^X$ and $v_1', v_2' \in \RR^{X'}$ we have
\begin{equation}
\label{kronip}
  \la v_1 \otimes v'_1, v_2 \otimes v'_2 \ra = \la v_1, v_2 \ra \la v'_1, v'_2 \ra,
\end{equation}
\vspace{-3.5mm}
\begin{equation}
\label{kron1a}
  (B_1 \otimes B_1')(v_1 \otimes v_1') = (B_1 v_1) \otimes (B_1'v_1').
\end{equation}


\section{The hypercubes and the matrices $\alpha_i$}
\label{sec:hyper}

In this section we recall the hypercube graph and some of its basic properties.
Fix a positive integer $D$. 
Let $X$ denote the Cartesian product 
$\{0,1\} \times \{0,1\} \times \cdots \times \{0,1\}$ ($D$ copies).
Thus $X$ is the set of sequences $(x_1, x_2, \ldots, x_D)$ such that $x_i \in \{0,1\}$ for $1 \le i \le D$.
For $x \in X$ and $1 \le i \le D$ let $x_i$ denote the $i$-th coordinate of $x$.
We interpret $\MX$ and the standard module $V=\RR^X$ as follows.
We abbreviate $\Mdva = \mat_{\{0,1\}}(\RR)$ and identify $\MX = (\Mdva)^{\otimes D}$.
We abbreviate $\RR^2 = \RR^{\{0,1\}}$ and identify $V=(\RR^2)^{\otimes D}$.
Let $Q_D$ denote the graph with vertex set $X$, and where two vertices are adjacent
if and only if they differ in exactly one coordinate. We call $Q_D$ the $D$-{\em cube} or the
$D$-dimensional {\em hypercube}. The graph $Q_D$ is connected and for $x,y \in X$ the distance $\partial(x,y)$ 
is the number of coordinates at which $x$ and $y$ differ. The diameter of $Q_D$ 
equals $D$. The graph $Q_D$ is bipartite with bipartition $X=X^+ \cup X^-$, where $X^+$ (resp.~$X^-$) 
is the set of vertices of $Q_D$ with an even (resp. odd) number of positive coordinates.
By \cite[p.~261]{BCN} $Q_D$ is distance-regular. 

\noindent
Let $\theta_0 > \cdots > \theta_D$ denote the eigenvalues of $Q_D$. By \cite[p.~261]{BCN} these 
eigenvalues and their multiplicities are given by
\begin{equation}
\label{eig_mult}
\theta_i = D-2i, \qquad \qquad \qquad \qquad m_i={D \choose i} \qquad \qquad \qquad \qquad
(0 \le i \le D).
\end{equation}

\begin{definition}
\label{def:i-adj}
{\rm For $1 \le i \le D$, vertices $x,y \in X$ are said to be $i$-{\em adjacent} whenever 
they differ in the $i$-th coordinate and are equal in all other coordinates. Define $\a_i \in \MX$ by 
$$
  (\a_i)_{x y} = \left\{ \begin{array}{ll}
                 1 & \hbox{if $x, y$ are $i$-adjacent}, \\
                 0 & \hbox{otherwise} \end{array} \right. \qquad (x,y \in X).
$$}
\end{definition}
From Definition \ref{def:i-adj} we routinely obtain the following result.

\begin{lemma}
\label{mat:lem1a}
With reference to Definition \ref{def:i-adj} the following {\rm (i)--(iii)} hold.
\begin{itemize}
\item[{\rm (i)}]
$\a_i \a_j = \a_j \a_i \; (1 \le i,j \le D)$;
\item[{\rm (ii)}]
$\a_i^2 = I \; (1 \le i \le D)$;
\item[{\rm (iii)}]
$A=\sum_{i=1}^D \a_i$.
\end{itemize}
\end{lemma}
We now describe the $\{\alpha_i\}_{i=1}^D$ from another point of view. 

\begin{definition}
\label{def:alpha}
{\rm 
Define $\a \in \Mdva$ by
$$
  \a = \left( \begin{array}{cc}
              0 & 1 \\
              1 & 0 
       \end{array} \right) .
$$
Observe that $\a^2=\mathnormal{1}$, where $\mathnormal{1}$ denotes the identity in $\Mdva$.}
\end{definition}

\begin{lemma}
\label{lem:alpha}
For $1 \le i \le D$ the matrix $\a_i$ from Definition \ref{def:i-adj} satisfies
\begin{equation}
\label{eq:alpha}
  \a_i = \mathnormal{1}^{\otimes (i-1)} \otimes \a \otimes \mathnormal{1}^{\otimes (D-i)}.
\end{equation}
\end{lemma}
\proof
Using the definition of tensor product in 
Section \ref{sec:kronecker}, along with Definition \ref{def:i-adj} and Definition \ref{def:alpha}, we find that 
for $x, y \in X$ the $(x,y)$-entry of the left-hand side of \eqref{eq:alpha} equals the $(x,y)$-entry of the 
right-hand side of \eqref{eq:alpha}. \qed 


\section{The matrices $\as_i$}
\label{sec:alpha_star}

We continue to discuss the hypercube $Q_D$ from Section \ref{sec:hyper}. In Section \ref{sec:hyper} 
we defined the matrices $\{ \a_i \}_{i=1}^D$. We now define some matrices 
$\{ \as_i \}_{i=1}^D$.

\begin{definition}
\label{def:i-adj_dual}
{\rm 
For $1 \le i \le D$ let $\as_i \in \MX$ denote the diagonal matrix with $(x,x)$-entry
$$
  (\as_i)_{xx} = \left\{ \begin{array}{ll}
                 1 & \hbox{if } \; x_i = 0, \\
                -1 & \hbox{if } \; x_i = 1 \end{array} \right. \qquad (x \in X).
$$}
\end{definition}
From Definition \ref{def:i-adj_dual} we routinely obtain the following result.

\begin{lemma}
\label{mat:lem1as}
With reference to Definition \ref{def:i-adj_dual} the following {\rm (i), (ii)} hold.
\begin{itemize}
\item[{\rm (i)}]
$\as_i \as_j = \as_j \as_i \; (1 \le i,j \le D)$;
\item[{\rm (ii)}]
$(\as_i)^2 = I \; (1 \le i \le D)$.
\end{itemize}
\end{lemma}

\smallskip \noindent
We now show that the matrices $\{\as_i\}_{i=1}^D$ satisfy an analog of Lemma \ref{lem:alpha}. 

\begin{definition}
\label{def:alpha_dual}
{\rm 
Define $\as \in \Mdva$ by 
$$
  \as = \left( \begin{array}{cc}
              1 & 0 \\
              0 & -1 
       \end{array} \right) .
$$
Observe that $(\as)^2=\mathnormal{1}$ and $\a \as = - \as \a$.}
\end{definition}
\begin{lemma}
\label{lem:i-adj_dual}
For $1 \le i \le D$ the matrix $\as_i$ from Definition \ref{def:i-adj_dual} satisfies
\begin{equation}
\label{eq:alpha_dual}
  \as_i = \mathnormal{1}^{\otimes (i-1)} \otimes \as \otimes \mathnormal{1}^{\otimes (D-i)}.
\end{equation}
\end{lemma}
\proof
Using the definition of tensor product in 
Section \ref{sec:kronecker}, along with Definition \ref{def:i-adj_dual} and Definition \ref{def:alpha_dual}, we find that 
for $x, y \in X$ the $(x,y)$-entry of the left-hand side of \eqref{eq:alpha_dual} equals the $(x,y)$-entry of the 
right-hand side of \eqref{eq:alpha_dual}. \qed 

\begin{lemma}
\label{mat:lem1}
With reference to Definitions \ref{def:i-adj} and \ref{def:i-adj_dual} the following {\rm (i), (ii)} hold.
\begin{itemize}
\item[{\rm (i)}]
$\a_i \as_j = \as_j \a_i$ if $i \ne j \; (1 \le i,j \le D)$;
\item[{\rm (ii)}]
$\a_i \as_i = - \as_i \a_i \; (1 \le i \le D)$.
\end{itemize}
\end{lemma}
\proof
Straightforward from Lemma \ref{lem:alpha} and Lemma \ref{lem:i-adj_dual}, using \eqref{kron1}
and $\a \as = -\as \a$. \qed


\section{An orthonormal $A$-eigenbasis for $V$}
\label{sec:basis}

We continue to discuss the hypercube $Q_D$ from Section \ref{sec:hyper}.
In this section we display an orthonormal basis for the standard module 
$V$ that consists of eigenvectors for $A$.

\begin{definition}
\label{def:uv}
{\rm Define $u,v \in \RR^2$ by 
$$
  u = {1 \over \sqrt{2}} \left( \begin{array}{c}
              1 \\
              1 
       \end{array} \right) , \qquad \qquad
  v = {1 \over \sqrt{2}} \left( \begin{array}{c}
              1 \\
              -1 
       \end{array} \right) .
$$
Observe that $u,v$ form an orthonormal basis for $\RR^2$.}
\end{definition}
We have a comment.

\begin{lemma}
\label{lem:uv}
The vectors $u,v$ from Definition \ref{def:uv} satisfy
$$
  \a u = u, \qquad \a v = -v, \qquad \as u = v, \qquad \as v = u,
$$
where $\a$ is from Definition \ref{def:alpha} and $\as$ is from Definition \ref{def:alpha_dual}.
\end{lemma}

\begin{definition}
\label{def:dual_vec}
{\rm 
For a subset $S \subseteq \{1,2, \ldots, D\}$ define $w_S \in V$ by
$$
  w_S = w_1 \otimes w_2 \otimes \cdots \otimes w_D,
$$
where
$$
  w_i = \left\{ \begin{array}{lll}
                 u & \hbox{if } \; i \not \in S,& \\
                 v & \hbox{if } \; i \in S & \end{array} \right. (1 \le i \le D).
$$
}
\end{definition}

\begin{lemma}
\label{lem:basis1}
The vectors 
\begin{equation}
\label{basis}
 w_S,  \qquad \qquad S \subseteq \{1,2, \ldots, D\}
\end{equation}
form an orthonormal basis for $V$.
\end{lemma}
\proof
The number of vectors in \eqref{basis} is $2^D$, and this number is the dimension of $V$.
Therefore it suffices to show that the vectors \eqref{basis} have square norm $1$ and are mutually orthogonal.
But this is the case by \eqref{kronip} and the observation below Definition \ref{def:uv}. \qed

\medskip \noindent
We now consider the actions of $\{\a_i\}_{i=1}^D, \: \{\as_i\}_{i=1}^D$ on the basis \eqref{basis}.

\begin{proposition}
\label{prop:basis2}
For $1 \le i \le D$ and $S \subseteq \{1,2, \ldots, D\}$ the action of $\a_i$  and $\as_i$ on $w_S$ is given by
$$
  \a_i w_S =   \left\{ \begin{array}{lll}
                 w_S & \hbox{if } \; i \not \in S,     & \\
                 -w_S & \hbox{if } \; i \in S, & \end{array} \right. 
\qquad \qquad
  \as_i w_S = \left\{ \begin{array}{lll}
                 w_{S \cup i} & \hbox{if } \; i \not \in S,     & \\
                 w_{S \setminus i} & \hbox{if } \; i \in S. & \end{array} \right.
$$ 
\end{proposition}
\proof
To compute $\a_i w_S$ use \eqref{kron1a}, Lemma \ref{lem:alpha}, Lemma \ref{lem:uv} and Definition \ref{def:dual_vec}.
To compute $\as_i w_S$ use \eqref{kron1a}, Lemma \ref{lem:i-adj_dual}, Lemma \ref{lem:uv} and Definition \ref{def:dual_vec}.
\qed

\begin{corollary}
\label{cor:basis3}
For $0 \le i \le D$ the vectors
\begin{equation}
\label{basis1}
 w_S, \qquad \qquad S \subseteq \{1,2, \ldots, D\}, \qquad \qquad |S|=i
\end{equation}
form a basis for $E_iV$.
\end{corollary}
\proof
In view of \eqref{de} and Lemma \ref{lem:basis1} it suffices to show that each vector $w_S$ from \eqref{basis1}
is contained in $E_iV$. Recall that $E_iV$ is the eigenspace of $A$ for the eigenvalue $\theta_i=D-2i$.
Using Lemma \ref{mat:lem1a}(iii) and Proposition \ref{prop:basis2},
$$
  A w_S = (D-|S|) w_S - |S| w_S = (D-2i) w_S.
$$ 
Therefore $w_S \in E_iV$ and the result follows. \qed


\section{A characterization of $\Al$}
\label{sec:char}

We continue to discuss the hypercube $Q_D$ from Section \ref{sec:hyper}.
For this graph we now give a characterization of $\Al$.

\begin{lemma}
\label{asal:lem0}
Pick distinct integers $i,j \; (1 \le i,j \le D)$ and $B \in \MX$.
Consider the following expression:
\begin{equation}
\label{C}
  \as_i \as_j B - \as_i B \as_j - \as_j B \as_i + B \as_i \as_j.
\end{equation}
Then for $x,y \in X$ the following {\rm (i), (ii)} hold.
\begin{itemize}
\item[{\rm (i)}]
The $(x,y)$-entry of \eqref{C} is equal to
\begin{equation}
\label{C1}
  ((\as_i)_{xx} - (\as_i)_{yy})((\as_j)_{xx} - (\as_j)_{yy}) B_{xy}.
\end{equation}
\item[{\rm (ii)}]
The $(x,y)$-entry of \eqref{C} is $0$ whenever $x=y$ or $x, y$ are adjacent.
\end{itemize}
\end{lemma}
\proof
(i) Use matrix multiplication, together with the fact that the matrices $\as_i$ and $\as_j$ are diagonal.

\smallskip \noindent
(ii) First assume $x=y$. Then the first two factors in \eqref{C1} are zero, so \eqref{C1} is zero. 
Next assume that $x,y$ are adjacent. Then there exists a unique integer $r \; (1 \le r \le D)$
such that $x,y$ are $r$-adjacent. By assumption $i \ne j$, so $r \ne i$ or $r \ne j$. 
If $r \ne i$ then the first factor in \eqref{C1} is zero. If $r \ne j$ then the second factor in \eqref{C1}
is zero. In any case \eqref{C1} is zero. \qed

\begin{lemma}
\label{asal:lem1}
For $B \in \MX$ the following {\rm (i), (ii)} are equivalent.
\begin{itemize}
\item[{\rm (i)}]
For all $x,y \in X$ that are not equal or adjacent, the $(x,y)$-entry of $B$ is zero.
\item[{\rm (ii)}]
For $1 \le i < j \le D$,
$$
  \as_i \as_j B - \as_i B \as_j - \as_j B \as_i + B \as_i \as_j = 0.
$$
\end{itemize}
\end{lemma}
\proof
(i) $\rightarrow$ (ii): Immediate from Lemma \ref{asal:lem0}. 

\medskip \noindent
(ii) $\rightarrow$ (i): 
Write $r=\partial(x,y)$ and note that $r \ge 2$. 
By construction $x,y$ differ in exactly $r$ coordinates. So there exist two distinct 
coordinates $i,j \; (i < j)$ at which $x,y$ differ. For these values of $i,j$ we apply Lemma \ref{asal:lem0}.
By assumption \eqref{C} is zero so \eqref{C1} is zero. But in \eqref{C1} the first two factors are nonzero
so the last factor $B_{xy}$ is zero. \qed

\begin{proposition}
\label{asal:prop1A}
For $B \in \MX$ the following {\rm (i), (ii)} are equivalent.
\begin{itemize}
\item[{\rm (i)}]
$B$ is $A$-like;
\item[{\rm (ii)}]
$B$ commutes with $A$ and 
$$
  \as_i \as_j B - \as_i B \as_j - \as_j B \as_i + B \as_i \as_j = 0 \qquad (1 \le i < j \le D).
$$
\end{itemize}
\end{proposition}
\proof
By Lemma \ref{asal:lem1} and the definition of an $A$-like matrix. \qed


\section{The symmetric $A$-like matrices for $Q_D$}
\label{sec:sal}

We continue to discuss the hypercube $Q_D$ from Section \ref{sec:hyper}.
For this graph we now describe the vector space $\Sal$.
We will give a basis for $\Sal$ and show that the dimension is $D+1$.
\begin{lemma}
\label{sal:lem1}
The following {\rm (i)--(iii)} hold.
\begin{itemize}
\item[{\rm (i)}]
$I \in \Sal$;
\item[{\rm (ii)}]
$\a_i \in \Sal$ for $1 \le i \le D$;
\item[{\rm (iii)}]
The matrices $I, \a_1, \a_2, \ldots, \a_D$ are linearly independent.
\end{itemize}
\end{lemma}
\proof
(i) This is clear.

\smallskip \noindent 
(ii) The matrix $\a_i$ is symmetric by Definition \ref{def:i-adj}, so it suffices to show that $\a_i \in \Al$.
By Lemma \ref{mat:lem1a}(i),(iii) we have $\a_i A = A \a_i$. For all $x,y \in X$ that are not equal or adjacent,
$x$ and $y$ are not $i$-adjacent so the $(x,y)$-entry of $\a_i$ is zero. We have shown $\a_i \in \Al$
and the result follows.

\smallskip \noindent
(iii) For each matrix in this list define the {\em support} to be the set of ordered pairs $(x,y)$ of 
vertices such that the $(x,y)$-entry is nonzero. These supports are nonempty and mutually disjoint.
Therefore the matrices are linearly independent. \qed

\begin{lemma}
\label{sal:lem2}
The following {\rm (i), (ii)} hold for $B \in \Sal$.
\begin{itemize}
\item[{\rm (i)}]
$B_{xx} = B_{yy}$ for all $x, y \in X$.
\item[{\rm (ii)}]
Pick $x,z \in X$ such that $\partial(x,z)=2$, and let $y,w$ denote the 
two vertices in $\G(x) \cap \G(z)$. 
Then $B_{xy}=B_{zw}$ and $B_{yz}=B_{wx}$.
\end{itemize}
\end{lemma}
\proof
(i) Since $Q_D$ is connected we may assume without loss that $x, y$ are adjacent. 
We have $BA=AB$ so $(BA)_{xy}=(AB)_{xy}$. By matrix multiplication,
$$
  (BA)_{xy} = \sum_{v \in X} B_{xv} A_{vy} = \sum_{v \in \G(y)} B_{xv}.
$$
By construction $x \in \G(y)$. For all $v \in \G(y) \backslash x$ we have
$\partial(x,v) = 2$, so $B_{xv}=0$. Therefore $(BA)_{xy}=B_{xx}$.
By a similar argument $(AB)_{xy}=B_{yy}$. The result follows.

\smallskip \noindent
(ii) Recall $BA=AB$ so $(BA)_{xz}=(AB)_{xz}$.
In this equation we expand each side using matrix multiplication and simplify the result using the fact that 
$B$ is in $\Sal$. This yields $B_{xy}+B_{wx} = B_{yz} + B_{zw}$.
In the above argument we replace $(x,y,z,w)$ by $(y,z,w,x)$ to obtain
$B_{yz}+B_{xy} = B_{zw} + B_{wx}$. 
Combining the above two equations we obtain $B_{xy}=B_{zw}$ and $B_{yz}=B_{wx}$. \qed

\begin{lemma}
\label{sal:lem2A}
Fix $x \in X$. Let $B$ denote a matrix in $\Sal$ such that $B_{xy}=0$ for all $y \in X$.
Then $B=0$.
\end{lemma}
\proof
For $z \in X$ we have $B_{zz}=0$ by Lemma \ref{sal:lem2}(i) and since $B_{xx}=0$.
We now show that $B_{zw}=0$ for all edges $zw$. We proceed as follows.
For the moment pick an edge $zw$. Since $Q_D$ is bipartite,
the distances $\partial(x,z)$ and $\partial(x,w)$ differ by $1$.
We claim that for all integers $m \; (1 \le m \le D)$,
$B_{zw}=0$ for all edges $zw$ such that $z \in \G_{m-1}(x)$ and $w \in \G_m(x)$.
To prove the claim we use induction on $m$.
First assume $m=1$. 
Then the claim holds by the assumptions of the lemma.
Next assume $m \ge 2$. Pick $v \in \G_{m-2}(x) \cap \G(z)$. Note that
$\partial(v,w)=2$; let $u$ denote the unique vertex in 
$\G(v) \cap \G(w)$ other than $z$. Applying Lemma \ref{sal:lem2}(ii)
to $v,z,w,u$ we find $B_{zw}=B_{vu}$. In this equation the right-hand side is 
zero by induction on $m$, so the left-hand side is zero, as desired.
The claim is proved and the result follows. \qed

\begin{corollary}
\label{sal:cor3}
The matrices
\begin{equation}
\label{basis-sal}
  I, \a_1, \a_2, \ldots, \a_D
\end{equation}
form a basis for $\Sal$. Moreover the dimension of $\Sal$ is $D+1$.
\end{corollary}
\proof
By Lemma \ref{sal:lem1} the matrices \eqref{basis-sal} are linearly independent and contained in $\Sal$.
It remains to show that the matrices \eqref{basis-sal} span $\Sal$. 
Pick $B \in \Sal$. We show that $B$ is in the span of \eqref{basis-sal}. 
Fix $x \in X$ and define real scalars $\varepsilon_0, \varepsilon_1, \ldots, \varepsilon_D$ as follows.
Define $\varepsilon_0 = B_{xx}$. For $1 \le i \le D$, define $\varepsilon_i = B_{xy}$, where 
$y=y_i$ is the unique vertex in $X$ that is $i$-adjacent to $x$. We show
\begin{equation}
\label{eq:sal} 
  B = \varepsilon_0 I + \sum_{i=1}^D \varepsilon_i \a_i.
\end{equation}
Let $C$ denote the left-hand side of \eqref{eq:sal} minus the right-hand side of \eqref{eq:sal},
and note that $C \in \Sal$.
By construction $C_{xy}=0$ for all $y \in X$, so $C=0$ in view of Lemma \ref{sal:lem2A}.
We have shown \eqref{eq:sal}. Therefore the matrices \eqref{basis-sal} span $\Sal$ and the result follows. \qed
 

\section{The antisymmetric $A$-like matrices for $Q_D$}
\label{sec:asal}

We continue to discuss the hypercube $Q_D$ from Section \ref{sec:hyper}.
For this graph we now describe the vector space $\Asal$.
We will give a basis for $\Asal$ and show that the dimension is ${D \choose 2}$.
We start with a comment. For $B \in \Al$ and $0 \le i \le D$, the space $E_iV$
is $B$-invariant since $B$ commutes with $A$.

\begin{lemma}
\label{asal:lem2}
For all $B \in \Asal$ we have $B E_0 V = 0$.
\end{lemma}
\proof
Note that $B E_0 V \subseteq E_0 V$.
Since $E_0 V$ has dimension $1$, there exists $\lambda \in \RR$ such that 
$(B - \lambda I)E_0 V =0$. We show that $\lambda = 0$. Pick a nonzero $v \in E_0 V$
and note that $\la v,v \ra \ne 0$. Since $B$ is antisymmetric, we get
$$
  \lambda \la v, v \ra = \la Bv, v \ra = -\la v, Bv \ra = -\lambda \la v, v \ra.
$$
Therefore $\lambda=0$ and the result follows. \qed

\medskip \noindent
To motivate our next result, pick $B \in \Asal$. Note that $E_1V$ is $B$-invariant. 
Consider the restriction $B|_{E_1 V}$. By Lemma \ref{lem:prelim2} this restriction is 
contained in the antisymmetric part of ${\rm End}(E_1V)$. 
Denote this part by ${\rm End}(E_1V)^{asym}$ and consider the restriction map
$\Asal \to {\rm End}(E_1V)^{asym}$, $B \mapsto B|_{E_1 V}$. We show that this restriction map
is an injection.

\begin{lemma}
\label{asal:lem3}
The restriction map 
$$
\begin{array}{ccl}
                 \Asal & \to & {\rm End}(E_1V)^{asym}\\
                 B & \mapsto & B|_{E_1 V} \end{array}
$$
is an injection.
\end{lemma}
\proof
Pick $B \in \Asal$ such that $BE_1V=0$.
We show that $B=0$.
We will do this in steps as follows.
We claim that $B E_\ell V = 0$ for $0 \le \ell \le D$.
Our proof is by induction on $\ell$. 
For $\ell=0$ the claim follows from Lemma \ref{asal:lem2},
and for $\ell=1$ the claim follows from our assumptions.
Next assume $\ell \ge 2$.
To show $B E_\ell V = 0$, by Corollary \ref{cor:basis3} it suffices to show $B w_S = 0$
for all subsets $S \subseteq \{1,2, \ldots, D\}$ such that $|S|=\ell$.
Let $S$ be given and pick distinct $i,j \in S$ with $i < j$. 
Define $P = S \setminus \{i,j\}$, $Q = S \setminus i$, $R = S \setminus j$.
Note that $|P|=\ell-2$, so $w_P \in E_{\ell-2}V$ by Corollary \ref{cor:basis3}.
Similarly $w_Q \in E_{\ell-1}V$ and $w_R \in E_{\ell-1}V$.
By these comments and the induction hypothesis, $B$ vanishes on each of 
$w_P$, $w_Q$, $w_R$.
By Proposition \ref{prop:basis2} we have $w_S = \as_i \as_j w_P$, $w_Q = \as_j w_P$, $w_R = \as_i w_P$. 
By Proposition \ref{asal:prop1A},
$$
  \as_i \as_j B - \as_i B \as_j - \as_j B \as_i + B \as_i \as_j = 0.
$$
In this equation we apply each side to $w_P$ and evaluate the result using the above comments to get $B w_S=0$.
We have shown $B E_\ell V = 0$ and the claim is proved. It follows that $B=0$. \qed

\medskip \noindent
We now show that $\as_i A \as_j - \as_j A \as_i \; (1 \le i < j \le D)$ form a basis
for $\Asal$. We start with a few comments about these expressions. To simplify the notation
we abbreviate 
\begin{equation}
\label{eq:Bij}
  B_{ij}=\as_i A \as_j - \as_j A \as_i \qquad \qquad (1 \le i < j \le D).
\end{equation}

\begin{lemma}
\label{asal:lem4a}
For $1 \le i < j \le D$,
$$
  B_{ij} = 2 \as_i \as_j(\a_i - \a_j).
$$
\end{lemma}
\proof
Routine using Lemma \ref{mat:lem1a}(iii), Lemma \ref{mat:lem1as} and Lemma \ref{mat:lem1}. \qed

\begin{lemma}
\label{asal:lem4b}
The following {\rm (i), (ii)} hold for $1 \le i < j \le D$ and $1 \le \ell \le D$.
\begin{itemize}
\item[{\rm (i)}]
Assume $\ell=i$ or $\ell=j$. Then $B_{ij} \a_\ell = - \a_\ell B_{ij}$.
\item[{\rm (ii)}]
Assume $\ell \ne i$ and $\ell \ne j$. Then $B_{ij} \a_\ell = \a_\ell B_{ij}$.
\end{itemize}
\end{lemma}
\proof
Use Lemma \ref{mat:lem1a}(i), Lemma \ref{mat:lem1}, and Lemma \ref{asal:lem4a}. \qed

\begin{lemma}
\label{asal:lem4c}
For $1 \le i < j \le D$ we have $B_{ij} \in \Asal$.
\end{lemma}
\proof
The matrix $B_{ij}$ is antisymmetric by \eqref{eq:Bij} and since each of $A$, $\as_i$ $\as_j$ is symmetric.
We show $B_{ij}$ commutes with $A$. 
Using Lemma \ref{mat:lem1a}(iii), Lemma \ref{asal:lem4a} and Lemma \ref{asal:lem4b} we find
$$
  B_{ij} A  - A B_{ij} = 2 B_{ij}(\a_i + \a_j) =  4 \as_i \as_j (\a_i^2 - \a_j^2) = 0,
$$
with the last equality holding since $\a^2_i=\a^2_j=I$. Therefore $B_{ij}$ commutes with $A$.

\smallskip \noindent
Pick $x,y \in X$ that are not equal or adjacent.
We show that the $(x,y)$-entry of $B_{ij}$ is zero.
By \eqref{eq:Bij} and since $\as_i$, $\as_j$ are diagonal,
$$
  (B_{ij})_{xy} = (\as_i)_{xx} A_{xy} (\as_j)_{yy} - (\as_j)_{xx} A_{xy} (\as_i)_{yy}.
$$
But $A_{xy}=0$ since $x,y$ are not adjacent, so $(B_{ij})_{xy} = 0$.
The result follows. \qed

\medskip \noindent
For $1 \le i < j \le D$ we now give the action of $B_{ij}$ on the basis \eqref{basis}. 

\begin{lemma}
\label{asal:lem5a}
For $1 \le i < j \le D$ and $S \subseteq \{1,2, \ldots, D\}$,
$$
  B_{ij} w_S = \left\{ \begin{array}{ll}
                           -4 w_{(S \cup j) \setminus i} & \hbox{if } \; i \in S \; \hbox{and } \; j \not \in S, \\
                           4w_{(S \cup i) \setminus j} & \hbox{if } \; i \not \in S \; \hbox{and } \; j \in S, \\
                           0 & \hbox{otherwise}.  \end{array} \right.
$$
\end{lemma}
\proof
Use Proposition \ref{prop:basis2} and Lemma \ref{asal:lem4a}. \qed

\begin{theorem}
\label{asal:thm5c}
The matrices 
\begin{equation}
\label{basis_asal}
   \as_i A \as_j - \as_j A \as_i, \qquad \qquad  1 \le i < j \le D
\end{equation}
form a basis for $\Asal$. Moreover the dimension of $\Asal$ is ${D \choose 2}$.
\end{theorem}
\proof
The number of elements in \eqref{basis_asal} is equal to ${D \choose 2}$.
The elements in \eqref{basis_asal} are contained in $\Asal$ by Lemma \ref{asal:lem4c}.
The elements in \eqref{basis_asal} are linearly independent; this can be verified using Lemma \ref{asal:lem5a} with $|S|=1$.
Therefore the dimension of $\Asal$ is at least ${D \choose 2}$.
To finish the proof it suffices to show that the dimension of $\Asal$ is at most ${D \choose 2}$.
By Lemma \ref{asal:lem3} the dimension of $\Asal$ is at most the dimension of ${\rm End}(E_1V)^{asym}$.
The dimension of ${\rm End}(E_1V)^{asym}$ is ${D \choose 2}$ since the dimension of $E_1V$ is $D$.
Therefore the dimension of $\Asal$ is at most ${D \choose 2}$. The result follows. \qed.

\begin{corollary}
\label{sal:cor4}
The following is a basis for $\Al$:
$$
  \{I, \a_1, \a_2, \ldots, \a_D\} \cup \{ \as_i A \as_j - \as_j A \as_i \mid 1 \le i < j \le D \}.
$$
Moreover the dimension of $\Al$ is $1 + D + {D \choose 2}$.
\end{corollary}
\proof
Recall that $\Al$ is a direct sum of $\Sal$ and $\Asal$. The result now follows from Corollary \ref{sal:cor3} and Theorem \ref{asal:thm5c}. \qed

\smallskip
\noindent
The following result might be of independent interest.

\begin{proposition}
\label{asal:cor6}
The restriction map 
$$
\begin{array}{ccl}
                 \Asal & \to & {\rm End}(E_1V)^{asym}\\
                 B & \mapsto & B|_{E_1 V} \end{array}
$$
is a bijection.
\end{proposition}
\proof
By Lemma \ref{asal:lem3} and since the dimensions of $\Asal$ and ${\rm End}(E_1V)^{asym}$ are equal. \qed


\end{document}